\documentclass[11pt]{article}
\usepackage[margin=1in]{geometry}
\usepackage{amsmath,amssymb,amsthm,bm}
\usepackage{booktabs}
\usepackage{enumitem}
\usepackage[numbers]{natbib}
\usepackage{xcolor}
\usepackage{tikz}
\usetikzlibrary{decorations.pathreplacing}
\definecolor{schedS}{HTML}{0072B2}
\definecolor{schedU}{HTML}{E69F00}
\definecolor{schedV}{HTML}{009E73}
\usepackage{algorithm}
\usepackage[noend]{algpseudocode}
\usepackage[colorlinks=true,linkcolor=blue,citecolor=blue,urlcolor=blue]{hyperref}

\newtheorem{theorem}{Theorem}[section]
\newtheorem{lemma}[theorem]{Lemma}
\newtheorem{proposition}[theorem]{Proposition}

\theoremstyle{definition}

\newtheorem{assumption}[theorem]{Assumption}
\newtheorem{example}[theorem]{Example}
\theoremstyle{remark}

\newcommand{\R}{\mathbb{R}}
\newcommand{\Z}{\mathbb{Z}}
\newcommand{\norm}[1]{\left\lVert #1\right\rVert}
\newcommand{\ip}[2]{\left\langle #1,#2\right\rangle}

\newcommand{\dist}{\operatorname{dist}}

\newcommand{\Fix}{\operatorname{Fix}}
\DeclareMathOperator*{\argmin}{arg\,min}

\newcommand{\xt}{\tilde{x}}
\newcommand{\xs}{x^{*}}
\newcommand{\Prop}{\mathcal{P}}

\title{PDFP: GPU-Accelerated Feasibility Pump via PDHG}
\author{Ruofeng Liu\thanks{MIT, Operations Research Center (alphakun@mit.edu).}
\and Haihao Lu\thanks{MIT, Sloan School of Management (haihao@mit.edu).}
\and Rahul Mazumder\thanks{MIT, Sloan School of Management (rahulmaz@mit.edu).}}
\date{}
\begin{document}
\maketitle

\begin{abstract}
The feasibility pump (FP) is among the most widely used primal heuristics for mixed-integer linear programming. It alternates between rounding an LP-feasible point and projecting the rounded point back onto the feasible region of the LP relaxation. Recent pump-based heuristics compute the projection with the primal--dual hybrid gradient method (PDHG), which is scalable on modern GPUs. They pose it as a linear program, which needs auxiliary variables and constraints for general integers. Once the projection runs on the GPU, the constraint-propagation rounding (CPR) of FP~2.0 takes a large share of each iteration. For the projection, we propose PDFP, a custom first-order method that solves the pump's projection problem directly, without auxiliary primal variables or constraints. We prove that it keeps the accelerated linear rate of restarted Halpern PDHG. For the rounding, we propose a new procedure that is more efficient than CPR and, under certain assumptions, still returns a CPR rounding. It also chooses rounded values by constraint violation, which can produce feasible roundings where nearest rounding cannot. Finally, our pump runs two independent trajectories concurrently on one GPU through CPU--GPU synchronization. This exploits the variance in trajectory length and lets the pump terminate sooner on average. Computational experiments on MIPLIB~2017 instances show early, promising results: these contributions speed up the discovery of the first incumbent solution. The code is publicly available at \url{https://github.com/MIT-Lu-Lab/PDFP}.
\end{abstract}

% \tableofcontents

%======================================================================
\section{Introduction}
\label{sec:intro}

Mixed-Integer Linear Programming (MILP) is a class of linear optimization problems in which some of the decision variables are integral:
\begin{equation}\label{eq:milp}
\min\{c^\top x \;:\; A_{\mathrm{ineq}}x\le b_{\mathrm{ineq}},\; A_{\mathrm{eq}}x=b_{\mathrm{eq}},\; \ell\le x\le u,\; x_j\in\Z\ (j\in I)\},
\end{equation}
where $I\subseteq\{1,\dots,n\}$ indexes the integer variables. Integer variables model indivisible quantities and logical decisions, such as whether to open a facility, commit a generator, or assign a crew to a flight. MILP is therefore widely used in scheduling, logistics, network design, production planning and energy systems~\citep{Wolsey2020}. It is also hard: solving a MILP is NP-hard, and even deciding whether a feasible solution exists is NP-complete~\citep{Karp1972}. In practice, MILPs are solved by branch-and-cut~\citep{LandDoig1960,PadbergRinaldi1991,Achterberg2007}. The solver enumerates a search tree over the integer variables and strengthens the linear programming (LP) relaxation of each node with cutting planes. It prunes every node whose LP bound is no better than the best feasible solution found so far, the \emph{incumbent}.

Primal heuristics, which search for good feasible solutions, are central to branch-and-cut~\citep{Berthold2006}. An early incumbent lets branch-and-cut prune more of its tree. On large or hard instances, where getting optimal solutions can be challenging, a good-quality feasible solution is very useful. Many instances of the MIPLIB~2017 library~\citep{MIPLIB2017} remain unsolved, and for some of them no feasible solution is known at all.

The focus of our paper is the feasibility pump (FP)~\citep{FGL2005}, which is among the most widely used primal heuristics. It keeps two points, an LP-feasible $\xs$ and an integral $\xt$, and alternates two steps. First, it rounds $\xs$ to obtain $\xt$. Then it projects $\xt$ back onto the feasible region of the LP relaxation, by minimizing the $\ell_1$ distance $\Delta(x,\xt)=\sum_{j\in I}|x_j-\xt_j|$. The loop stops when the two points meet, and random perturbations break the cycles it can enter. Later variants extend the pump to general integers with a two-stage scheme~\citep{BFL2007} and add the objective to the distance~\citep{AchterbergBerthold2007}. FP~2.0~\citep{FP20} replaced nearest rounding with constraint-propagation rounding (CPR), a diving-like procedure that fixes one variable at a time and propagates each fixing through the constraints. FP~2.0 has been the reference variant of the pump since. In a classical implementation, every projection is a simplex solve, which dominates the time of an iteration. The rounding work of CPR is small by comparison.

In recent years, the primal--dual hybrid gradient method (PDHG)~\citep{CP} has been gaining popularity as an algorithm for solving LP~\citep{PDLP2021,ApplegateHinderLuLubin2023,cuPDLP2023,HalpernPDHG2024,cuPDLPx2025}. It is a matrix-free first-order method, meaning that it does not need to solve any linear system, and its computational bottleneck is matrix--vector multiplication. It is also memory efficient, as it does not need to store the factorization of any matrix. PDLP~\citep{PDLP2021} made restarted PDHG a practical LP method. Using the parallelization capabilities of matrix--vector multiplication, this method can be scalable on modern GPUs, which the GPU solvers cuPDLP~\citep{cuPDLP2023,lu2023cupdlp}, cuPDLPx~\citep{HalpernPDHG2024,cuPDLPx2025} and D-PDLP~\cite{li2026d} exploit.

Recently, these new LP algorithms are being used to speed up MILP heuristics. In particular, Scylla~\citep{Scylla2023} uses PDHG for the projections of an FP-like loop on the CPU. NVIDIA's cuOpt~\citep{Corduk2025cuOpt} uses it for the projections of a GPU feasibility pump whose bound propagation and bulk rounding also run on the GPU. Some works use it for the root relaxation alone~\citep{Kempke2026LowPrecision,CHAP2026}. With a PDHG projector, the cost of a projection grows with the accuracy requested. It often suffices to do inexact projections. This allows us to spend the time saved on more pump iterations. This strategy is used by Scylla, with a tolerance that starts loose and tightens geometrically over the run~\citep{Scylla2023}.

A pump with a PDHG projector still faces two issues. First, prior pumps, Scylla and cuOpt included, solve the projection as an LP in which the absolute values of the $\ell_1$ distance are linearized. For general integers, this LP needs auxiliary variables and constraints, which enlarge the operator that every PDHG iteration applies. Second, CPR is sequential, with one propagation per decision, and with the projection on the GPU the rounding takes a large share of each pump iteration. Such a pump also offers two opportunities. A PDHG projection runs many iterations, and an intermediate iterate may already round to a feasible point. When the rounding runs on the CPU, the GPU is idle during every rounding.

Here we advance the frontiers of GPU-powered algorithms for FP. We propose and implement an enhanced feasibility pump that tackles both issues and uses both opportunities. Our aim is a pump that finds feasible solutions on more instances, and sooner. Our contributions are as follows.
\begin{enumerate}[nosep]
  \item \textbf{PDFP, a custom first-order solver for the pump's projection problem, and its convergence} (Sections~\ref{sec:pdfp} and~\ref{sec:theory}). We design PDFP, which dualizes the $\ell_1$ distance term of the projection and runs restarted Halpern PDHG with reflection (r$^2$HPDHG)~\citep{HalpernPDHG2024} on the resulting saddle-point problem. PDFP solves the original projection problem without auxiliary primal variables or additional constraints. So every iteration works with the original constraint matrix, also on instances with general integers. We prove that PDFP converges at the accelerated linear rate that r$^2$HPDHG attains on LP, so this guarantee carries over to the pump's projection problem.
  \item \textbf{Guarded two-phase CPR, a rounder that fixes many variables at once} (Section~\ref{sec:cpr}). Our rounder fixes, in a single propagation, all integer coordinates whose LP value is already integral, and falls back to sequential CPR when this batch contradicts itself. We prove that, under certain assumptions, every output of our rounder is a sequential CPR rounding, unlike that of bulk rounding~\citep{Corduk2025cuOpt}. So the batch speeds up the rounder, and the rounder still returns a CPR rounding. We also design an online cost rule that learns, per run, whether trying the batch pays, which produces additional speedups. Finally, CPR uses the constraints to shape each variable's domain but never to choose its value. We propose a rule that rounds each variable to the value in its domain that adds the least constraint violation, given the variables fixed so far. We observed that this rule can produce feasible roundings when the nearest-value rule of FP~2.0 cannot.
  \item \textbf{Further enhancements to the GPU pump} (Section~\ref{sec:enhancements}). We stop the pump as soon as a PDHG iterate rounds to a feasible point. This allows us to capture mid-projection feasible points that we would miss otherwise. We also run two trajectories of the pump concurrently, with independent seeds, on one GPU through CPU--GPU synchronization: the projections of each tend to fall in the rounding phase of the other. Where the rounding dominates, the second trajectory therefore comes at nearly zero cost in wall-clock time. This exploits the variance in trajectory length and allows our pump to terminate sooner on average. Finally, the pump detects cycles by hashing the visited points, which saves memory and time.
\end{enumerate}
We evaluate the pump by its time to a first feasible solution on MIPLIB~2017 instances with at least one million nonzeros. We compare it with the FP~2.0 reference implementation of Salvagnin and Mexi~\citep{FeasPumpCollection}, run with cuPDLPx as its LP solver, and with NVIDIA's cuOpt~\citep{Corduk2025cuOpt}. Our pump finds feasible solutions on more instances than both, and finds them faster: 1.18$\times$ faster than the reference implementation and 1.72$\times$ faster than cuOpt in shifted geometric mean.

%We also report, with mechanisms, the enhancements that did not work (Section~\ref{sec:negative}).

The rest of the paper is organized as follows. Section~\ref{sec:related} reviews related work. Section~\ref{sec:background} states the feasibility pump loop, FP~2.0's constraint-propagation rounding and restarted Halpern PDHG. Section~\ref{sec:pdfp} presents PDFP, and Section~\ref{sec:theory} proves its convergence. Section~\ref{sec:cpr} presents our rounder, guarded two-phase CPR with violation-guided rounding. Section~\ref{sec:enhancements} describes the further enhancements of our pump, and Section~\ref{sec:results} reports the computational results. Section~\ref{sec:conclusion} concludes.

%======================================================================
\section{Related work}
\label{sec:related}

In this section, we review the work that our pump builds on and the methods it is compared with. Our designs draw on four lines of research: the feasibility pump and its variants, propagation-based rounding, first-order LP methods and their convergence theory, and the use of first-order LP solvers inside primal heuristics. We close with hybrid and LP-free heuristics, which pursue the same goal, a first feasible solution, by other means.

\paragraph{The feasibility pump lineage.}
FP was introduced for binary MILPs by \citet{FGL2005} and extended to general integers by \citet{BFL2007}, whose three-stage scheme pumps on the binaries first, then on all integer variables, and ends with an enumeration stage. \citet{AchterbergBerthold2007} made the distance objective-aware. \citet{FP20} replaced rounding by constraint-propagation rounding, which is the baseline throughout this paper. \citet{BertholdLodiSalvagnin2019} survey the first decade. Several lines revisit individual ingredients. \citet{DeSantis2013,DeSantis2014} interpret FP as a Frank--Wolfe method on a concave merit function and replace the $\ell_1$ distance with other integrality measures. \citet{Geissler2017} view FP as a penalty alternating direction method. \citet{DeyEtAl2018} analyze and improve the randomization step. \citet{Boland2014} improve the rounding step by choosing more carefully the point that is rounded. \citet{BertholdMexiSalvagnin2023} combine several reference vectors in the distance. A reference implementation of FP~2.0, FeasPumpCollection~\citep{FeasPumpCollection}, is public. Kernel Pump~\citep{KernelPump2026} runs FP~2.0 on a kernel of binary variables chosen from the LP relaxation; since it changes neither the rounder nor the projector, our pump can serve as its inner FP~2.0.

\paragraph{Propagation-based rounding and fix-and-propagate.}
CPR~\citep{FP20} rounds one variable at a time and propagates each fixing with bound strengthening, scheduled by a propagation engine over row propagators~\citep{Apt2003}. Related CPU heuristics include shift-and-propagate~\citep{BertholdHendel2015} and fix-propagate-repair~\citep{Salvagnin2024FPR}, which runs LP-free dives with a WalkSAT-like repair. \citet{Sofranac2022} parallelize domain propagation over sparse matrices on the GPU. The closest prior art to our batched rounder is NVIDIA cuOpt's GPU fix-and-propagate~\citep{Corduk2025cuOpt}, which rounds bulks of variables together, propagates once per bulk and backtracks when a bulk fails. Its output differs from that of sequential rounding by design. Our batch is the set of LP-integral coordinates, and we characterize exactly what the guarded rounder computes relative to sequential CPR (Section~\ref{sec:cpr-exact}). Our rounder runs on the CPU, but the batch is an algorithmic change that does not depend on the device and could be combined with GPU propagation.

\paragraph{PDLP.}
PDLP~\citep{PDLP2021} demonstrated the potential of restarted PDHG~\citep{CP} as a practical LP method. GPU implementations such as cuPDLP~\citep{cuPDLP2023,lu2023cupdlp} exploit the parallelism of sparse matrix--vector products and have further demonstrated its practical performance. These advances have attracted interest across the optimization solver community, with variants incorporated into NVIDIA cuOpt, Gurobi, COPT, HiGHS, FICO Xpress and Knitro. Subsequent developments have further improved performance and scalability: cuPDLPx~\citep{HalpernPDHG2024,cuPDLPx2025} implements restarted Halpern PDHG (rHPDHG), while D-PDLP~\citep{li2026d} provides a distributed implementation across multiple GPUs.

On the theory side, PDHG can be interpreted as a preconditioned proximal point method, and its ergodic iterates converge at the sublinear rate $O(1/k)$ for convex--concave saddle-point problems~\citep{HeYuan2012,LuYang2023Unified}. Under metric subregularity, which holds for LP on bounded sets, the last iterate of PDHG converges linearly, requiring $O(\kappa^2\log(1/\varepsilon))$ iterations, where $\kappa$ is the ratio of a smoothness constant to a sharpness constant of the LP~\citep{LuYang2022IDS}. Restarts improve this complexity to $O(\kappa\log(1/\varepsilon))$, matching a lower bound for a wide class of primal--dual methods. Without restarts, the quadratic dependence on $\kappa$ cannot be improved in general~\citep{ApplegateHinderLuLubin2023}. Restarted Halpern PDHG attains the same accelerated rate $O(\kappa\log(1/\varepsilon))$, with a further factor-of-two improvement from reflection~\citep{HalpernPDHG2024}. Our convergence result (Section~\ref{sec:theory}) extends this analysis to the pump's $\ell_1$ saddle-point problem by establishing the required sharpness property using the upper-Lipschitz property of polyhedral multifunctions~\citep{Robinson1981}.

\paragraph{First-order projections inside heuristics.}
Scylla~\citep{Scylla2023} is the closest prior work on the projection side: an FP-like fix-propagate-and-project loop with warm-started PDLP projections whose tolerance starts loose and tightens geometrically. It runs on the CPU, and its authors leave the parallel potential of PDHG as future work. The pump of NVIDIA's cuOpt~\citep{Corduk2025cuOpt} warm-starts PDLP on the GPU, caps each projection in time and continues from the point reached at the cap. CHAP~\citep{CHAP2026} solves only the root LP on the GPU, while its pump runs on the CPU. \citet{Kempke2026LowPrecision} guide fix-and-propagate with one low-precision PDLP root solve, and \citet{Xiong2026SCUC} use a GPU implementation of a first-order LP method for successive fixing on unit-commitment models. All of these solve LPs. PDFP instead solves the $\ell_1$ projection itself rather than its LP reformulation (Section~\ref{sec:pdfp}).

\paragraph{Hybrid and LP-free heuristics.}
Feasibility Jump~\citep{LutebergetSartor2023} is an LP-free weighted local search, and cuOpt pairs its pump with GPU local search. CHAP~\citep{CHAP2026} coordinates GPU and CPU heuristics through a shared solution pool. Our use of idle device time (Section~\ref{sec:racing}) is finer-grained: it draws a second trajectory of the same pump on the same device.

%======================================================================
\section{Background}
\label{sec:background}

In this section, we fix the notation and review the two methods our work builds on. Section~\ref{sec:fp20} states the FP~2.0 loop on which our pump builds, together with its constraint-propagation rounding, whose number of propagations our rounder reduces. Section~\ref{sec:pdhg} reviews restarted Halpern PDHG, on which both of our projectors are built.

\paragraph{Notation.} We use the notation of~\eqref{eq:milp}: $A\in\R^{m\times n}$ stacks $A_{\mathrm{ineq}}$ and $A_{\mathrm{eq}}$, the bounds $\ell,u$ may be infinite, $X=[\ell,u]$ is the box, and $\Omega=\{x\in X: A_{\mathrm{ineq}}x\le b_{\mathrm{ineq}},\, A_{\mathrm{eq}}x=b_{\mathrm{eq}}\}$ is the feasible region of the LP relaxation. The integer set $I$ splits into binaries $\mathcal{B}$ (integer columns with $u_j-\ell_j=1$) and general integers $\mathcal{G}$. A point is \emph{integral} if $x_j\in\Z$ for all $j\in I$, and for $J\subseteq\{1,\dots,n\}$, $x_J:=(x_j)_{j\in J}$ denotes the subvector of $x$ indexed by $J$. For a closed convex set $C$, $\Pi_C$ denotes the Euclidean projection onto $C$ and $N_C(x)$ the normal cone of $C$ at $x\in C$. We write $(\cdot)_+=\max\{\cdot,0\}$ componentwise, $\Fix(T):=\{z:T(z)=z\}$ for the set of fixed points of an operator $T$, and $\lambda_{\max}(\cdot)$ for the largest eigenvalue of a symmetric matrix. We use $\norm{\cdot}_2$ for the Euclidean norm of a vector and the spectral norm of a matrix, and drop the subscript where no confusion can arise. For a symmetric positive definite matrix $P$, $\ip{u}{v}_P:=\ip{u}{Pv}$, $\norm{u}_P:=\ip{u}{u}_P^{1/2}$ and $\dist_P(z,C):=\min_{c\in C}\norm{z-c}_P$; $\dist_2$ is the Euclidean distance.

\subsection{The feasibility pump and FP~2.0}
\label{sec:fp20}

Algorithm~\ref{alg:fp} states the FP~2.0 loop~\citep{FGL2005,FP20} on which our pump builds. The projection step minimizes $\Delta(\cdot,\xt)$ over $\Omega$, and implementations of the pump solve it as a linear program by linearizing the absolute values~\citep[Sect.~2]{FP20}. When $\xt_j$ sits at a bound of $[\ell_j,u_j]$, the term $|x_j-\xt_j|$ is linear on $X$ ($x_j-\ell_j$ or $u_j-x_j$). Every other integer coordinate needs an auxiliary variable $d_j$ and two rows, which gives the LP
\begin{equation}\label{eq:proj-lp}
\min_{x\in\Omega,\ d}\ \sum_{j\in I:\,\xt_j=\ell_j}(x_j-\ell_j)+\sum_{j\in I:\,\xt_j=u_j}(u_j-x_j)+\sum_{j\in I:\,\ell_j<\xt_j<u_j}d_j
\quad\text{s.t.}\quad d_j\ge x_j-\xt_j,\ \ d_j\ge\xt_j-x_j .
\end{equation}
On pure-binary instances every $\xt_j$ sits at a bound, so~\eqref{eq:proj-lp} has no auxiliary variables and its objective coefficients are $\pm1$. General integers are also harder for the pump than binaries, because rounding them must choose a value and not only a direction~\citep[Sect.~2.1]{FP20}. The scheme of Bertacco et al.~\citep{BFL2007} therefore pumps first on the binaries, with the integrality of the general integers relaxed, and then on the full model; we call these two stages the two-stage scheme.

\begin{algorithm}[t]
\caption{The feasibility pump loop of FP~2.0.}
\label{alg:fp}
\begin{algorithmic}[1]
\Require MILP~\eqref{eq:milp}; projector \textsc{Proj}; rounder \textsc{Round}; flip parameter $T$; restart period $R$; time limit
\State $\xs\gets$ an optimal solution of the LP relaxation (original objective)
\If{$\xs$ is integral and feasible} \Return $\xs$ \EndIf
\State $\xt\gets\textsc{Round}(\xs,0)$;\quad $\mathcal{V}\gets\emptyset$ \Comment{$\mathcal{V}$: visited integer points}
\For{$k=1,2,\dots$ until the time limit}
  \If{$\xt\in\Omega$} \Return $\xt$ \EndIf
  \State $\xs\gets\textsc{Proj}(\xt)\in\argmin\{\Delta(x,\xt):x\in\Omega\}$ \Comment{warm-started}
  \If{$\xs$ is integral and feasible} \Return $\xs$ \EndIf
  \State $\hat x\gets\textsc{Round}(\xs,k)$
  \If{$\hat x=\xt$} flip the $TT\sim\mathcal{U}\{\tfrac T2,\dots,\tfrac{3T}2\}$ integer coordinates with largest $|\xs_j-\xt_j|$ \EndIf
  \If{$\hat x\in\mathcal{V}$ \textbf{or} $R$ iterations since the last perturbation}
     \State perturb $\hat x$; repeat while $\hat x\in\mathcal{V}$ (bounded number of times)
  \EndIf
  \State $\mathcal{V}\gets\mathcal{V}\cup\{\xt\}$;\quad $\xt\gets\hat x$
\EndFor
\end{algorithmic}
\end{algorithm}

The perturbation flips binaries at random as in \citet{FGL2005}, and reassigns general integers at random over their domains~\citep[Sect.~5.1]{FP20}. As in the reference implementation~\citep{FeasPumpCollection}, the pump remembers every visited integer point and perturbs again while the perturbed point has been visited.

\paragraph{Constraint-propagation rounding.}
CPR~\citep[Sect.~4 and~6]{FP20} ranks the integer variables by increasing fractionality $|\xs_j-[\xs_j]|$, binaries before general integers. Following that order, it fixes each variable that is not yet fixed to the value of $\xs_j$ rounded within the variable's current domain, and propagates the fixing through the constraints by bound strengthening until no bound changes. We call this rounding of a variable a \emph{decision}. A variable that earlier propagation has already fixed is not decided; its value is read off its singleton domain. At the end every domain is a singleton and $\xt$ is read off. To reduce cycling, FP~2.0 randomly swaps $10\%$ of the variables in the ranking after the tenth pump iteration~\citep[Sect.~6]{FP20}. Algorithm~\ref{alg:seqcpr} states the procedure, which we call \emph{sequential CPR}. In it, $\textsc{Propagate}(B)$ applies bound strengthening to the rows until no bound of the box $B$ changes~\citep[Fig.~3]{FP20}. It returns the new box together with an \emph{infeasibility flag}, which it raises when some row can no longer be satisfied. Sequential CPR ignores the flag and always returns a rounding, which may be infeasible; the guarded rounder of Section~\ref{sec:cpr} reads it. The continuous coordinates of $\xt$ keep their LP values.

\begin{algorithm}[t]
\caption{Sequential CPR (FP~2.0).}
\label{alg:seqcpr}
\begin{algorithmic}[1]
\Require LP point $\xs$; pump iteration $k$; box $X=[\ell,u]$
\State $\pi\gets$ the integer variables sorted by fractionality $|\xs_j-[\xs_j]|$, binaries first
\If{$k\ge10$} apply $\lfloor|I|/10\rfloor$ random transpositions to $\pi$ \EndIf
\State $B\gets X$ with $B_j=[\lceil\ell_j\rceil,\lfloor u_j\rfloor]$ for $j\in I$;\quad $\xt\gets\xs$
\For{each $j$ in the order $\pi$}
  \If{$B_j$ is a singleton $\{v\}$} $\xt_j\gets v$ \Comment{fixed by earlier propagation}
  \Else
    \State $v\gets$ the integer in $B_j$ nearest to $\xs_j$ \Comment{the plain decision rule}
    \State $B\gets\textsc{Propagate}(B\cap\{x_j=v\})$;\quad $\xt_j\gets v$
  \EndIf
\EndFor
\State \Return $\xt$
\end{algorithmic}
\end{algorithm}

Sequential CPR runs one propagation per decision, and these propagations run one after another, because each decision reads the box that the previous propagation left. The work inside each propagation can be parallelized, on the GPU as well~\citep{Sofranac2022,Corduk2025cuOpt}, but the number of propagations falls only if several variables are rounded at once, as in cuOpt's bulk rounding and in our batch (Section~\ref{sec:cpr}).

\subsection{PDHG and the restarted reflected Halpern iteration}
\label{sec:pdhg}

PDHG~\citep{CP} is the base algorithm of PDLP~\citep{PDLP2021}, cuPDLP~\citep{cuPDLP2023} and cuPDLPx~\citep{cuPDLPx2025}. We review it on the LP relaxation of~\eqref{eq:milp}, written as the saddle-point problem
\begin{equation}\label{eq:lp-saddle}
\min_{x\in X}\ \max_{y\in Y}\ L(x,y):=c^\top x-\ip{y}{Ax-b},\qquad Y:=\{y\in\R^m:\ y_i\le0\ \text{for the rows of }A_{\mathrm{ineq}}\}.
\end{equation}
Maximizing $L$ over $y\in Y$ returns $c^\top x$ if $x$ satisfies the constraints and $+\infty$ otherwise. One PDHG update $z^{k+1}=T(z^k)$, where $z^k:=(x^k,y^k)$, is given by
\begin{equation}\label{eq:pdhg}
x^{k+1}=\Pi_X\big(x^k-\eta(c-A^\top y^k)\big),\qquad
y^{k+1}=\Pi_Y\big(y^k-\sigma(A(2x^{k+1}-x^k)-b)\big),
\end{equation}
where $\eta,\sigma>0$ are the primal and dual step sizes. Each update multiplies one vector by $A^\top$ (in the $x$-update) and one vector by $A$ (in the $y$-update), plus vector operations. The method is matrix-free, and since sparse matrix--vector products parallelize well, it is parallelizable and scalable on modern GPUs.

The optimality conditions of~\eqref{eq:lp-saddle} are the monotone inclusion $0\in F(z)$ with $F(x,y):=\big(c-A^\top y+N_X(x),\ Ax-b+N_Y(y)\big)$. Following the proximal-point interpretation of PDHG~\citep{HeYuan2012,LuYang2023Unified}, the update $z^{k+1}=T(z^k)$ is the unique solution of
\begin{equation}\label{eq:ppm}
0\in F(z^{k+1})+P(z^{k+1}-z^k),\qquad P:=\begin{pmatrix}\eta^{-1}I & A^\top\\ A & \sigma^{-1}I\end{pmatrix}.
\end{equation}
When $\eta\sigma\norm{A}_2^2<1$, $P$ is positive definite and $T$ is firmly nonexpansive in the $P$-norm. The fixed points of $T$ are exactly the saddle points of~\eqref{eq:lp-saddle}. We call $\norm{z-T(z)}_P$ the \emph{fixed-point residual}.

\paragraph{Restarted reflected Halpern iteration.} Since $T$ is firmly nonexpansive, its reflection $2T-I$ is nonexpansive in the $P$-norm and has the same fixed points as $T$~\citep[Prop.~4.4]{BauschkeCombettes2017}. The restarted reflected Halpern PDHG (r$^2$HPDHG) of \citet{HalpernPDHG2024} applies the Halpern iteration~\citep{Halpern1967} to $2T-I$ and resets its anchor at restarts. With outer counter $n$ and inner counter $t$, it generates
\begin{equation}\label{eq:halpern}
z^{n,t+1}=\frac{t+1}{t+2}\big(2T(z^{n,t})-z^{n,t}\big)+\frac{1}{t+2}\,z^{n,0},
\end{equation}
and, if a restart is triggered after inner iteration $t_n$, it starts the next outer iteration from $z^{n+1,0}:=T(z^{n,t_n})$. A restart is triggered when the fixed-point residual has decayed enough. On feasible and bounded LPs, this scheme converges linearly, at an accelerated rate~\citep{HalpernPDHG2024}.

The cuPDLPx solver implements this scheme with a constant step size, diagonal preconditioning~\citep{Ruiz2001,PockChambolle2011}, a PID controller that balances the primal and dual step sizes at restarts, and termination on relative primal and dual residuals and duality gap~\citep{cuPDLPx2025}. Our PDLP projector is a Julia port of cuPDLPx, and PDFP (Section~\ref{sec:pdfp}) reuses its machinery. For the pump, one feature of these methods matters most: the cost of a solve grows with the accuracy requested, so the tolerance becomes a parameter the pump can set for each projection (Section~\ref{sec:tolerance}).

%======================================================================
\section{PDFP: a first-order projector for the pump's \texorpdfstring{$\ell_1$}{l1} problem}
\label{sec:pdfp}

In this section, we present PDFP, the projector of our pump. The projection step of the pump minimizes an $\ell_1$ distance over the feasible region of the LP relaxation, and prior pumps solve it as a linear program in which the absolute values are linearized. PDFP instead solves the $\ell_1$ problem itself: it dualizes the distance term and runs r$^2$HPDHG on the resulting saddle-point problem. This matters on instances with general integers, where the linear program needs an auxiliary variable and two rows for every integer coordinate whose target lies strictly between its bounds, while PDFP needs no auxiliary primal variables or additional constraints. Section~\ref{sec:pdfp-iter} derives the saddle-point problem and the algorithm, Section~\ref{sec:pdfp-precond} explains how preconditioning keeps the distance exact, and Section~\ref{sec:pdfp-term} gives the termination criteria.

The pump's projection is
\begin{equation}\label{eq:proj}
\min_{x\in\Omega}\ \varphi(x),\qquad \varphi(x):=\sum_{j=1}^n q_j\,|x_j-\xt_j|,
\end{equation}
where $q_j=1$ for $j\in I$ and $q_j=0$ otherwise; preconditioning changes these weights (Section~\ref{sec:pdfp-precond}). Prior implementations of the pump, including those whose projector is PDHG~\citep{Scylla2023,Corduk2025cuOpt}, solve~\eqref{eq:proj} as the linear program~\eqref{eq:proj-lp}, whose auxiliary variables and rows enlarge the matrix that every PDHG iteration multiplies.

\subsection{The saddle-point problem and the algorithm}
\label{sec:pdfp-iter}

Writing $q_j|x_j-\xt_j|=\max_{|s_j|\le q_j}s_j(\xt_j-x_j)$, we obtain~\eqref{eq:proj} as the primal problem of the saddle-point problem
\begin{equation}\label{eq:pdfp}
\min_{x\in X}\ \max_{s\in Q,\ y\in Y}\ L(x,s,y):=\ip{s}{\xt-x}-\ip{y}{Ax-b},\qquad Q:=[-q,q],
\end{equation}
with $Y$ as in~\eqref{eq:lp-saddle}. Maximizing $L$ over $s\in Q$ returns $\varphi(x)$, and maximizing it over $y\in Y$ returns $0$ if $x\in\Omega$ and $+\infty$ otherwise. Minimizing $L$ over $x\in X$ instead gives the dual problem
\begin{equation}\label{eq:pdfp-dual}
\max_{s\in Q,\ y\in Y}\ \ip{\xt}{s}+\ip{b}{y}-p\big(s+A^\top y;\ell,u\big),
\end{equation}
where $p(r;\ell,u):=\sup_{x\in X}r^\top x$ is the support function of the box $X$, possibly $+\infty$. Unlike~\eqref{eq:proj-lp}, problem~\eqref{eq:pdfp} keeps every term $q_j|x_j-\xt_j|$ as it is, whether $\xt_j$ lies at a bound of $x_j$ or strictly between its bounds. It needs no auxiliary primal variables or additional constraints: its only new variables are the multipliers $s$, which are updated by a projection onto a box.

Problem~\eqref{eq:pdfp} has the form~\eqref{eq:lp-saddle} with $c=0$, the matrix $A$ replaced by the stacked matrix $K:=\bigl(\begin{smallmatrix}I\\ A\end{smallmatrix}\bigr)$, the right-hand side $b$ replaced by $(\xt,b)$, and the dual set $Y$ replaced by $Q\times Y$. PDHG for~\eqref{eq:pdfp}, $z^{k+1}=T(z^k)$ with $z^k:=(x^k,s^k,y^k)$, therefore reads
\begin{equation}\label{eq:pdfp-step}
\begin{aligned}
x^{k+1}&=\Pi_X\big(x^k+\eta(s^k+A^\top y^k)\big),\\
s^{k+1}&=\Pi_Q\big(s^k+\sigma(\xt-2x^{k+1}+x^k)\big),\\
y^{k+1}&=\Pi_Y\big(y^k-\sigma(A(2x^{k+1}-x^k)-b)\big).
\end{aligned}
\end{equation}
Everything in Section~\ref{sec:pdhg} carries over with $A$ replaced by $K$: the update is the proximal-point step~\eqref{eq:ppm} with $P:=\begin{pmatrix}\eta^{-1}I&K^\top\\K&\sigma^{-1}I\end{pmatrix}$, and it requires $\eta\sigma\norm{K}_2^2<1$, where $\norm{K}_2^2=1+\norm{A}_2^2$. Each iteration still multiplies one vector by $A^\top$ and one vector by $A$; the $s$-update adds $O(n)$ work. PDFP runs r$^2$HPDHG on this operator $T$ (Algorithm~\ref{alg:pdfp}). It restarts when the fixed-point residual has decayed by a fixed factor $\beta\in(0,1)$,
\begin{equation}\label{eq:restart}
\norm{z^{n,t}-T(z^{n,t})}_P\le\beta\,\norm{z^{n,0}-T(z^{n,0})}_P,
\end{equation}
and it terminates at the first PDHG point $T(z^{n,t})$ that meets the criteria of Section~\ref{sec:pdfp-term}.

\begin{algorithm}[t]
\caption{PDFP: r$^2$HPDHG for the projection~\eqref{eq:pdfp}.}
\label{alg:pdfp}
\begin{algorithmic}[1]
\Require target $\xt$; initial point $z^{0,0}$; tolerances $\varepsilon_g,\varepsilon_r$; restart parameter $\beta\in(0,1)$
\State $n\gets0$
\Repeat
  \State $t\gets0$
  \Repeat
    \State $z^{n,t+1}\gets\frac{t+1}{t+2}\big(2T(z^{n,t})-z^{n,t}\big)+\frac{1}{t+2}\,z^{n,0}$;\quad $t\gets t+1$ \Comment{$T$ is~\eqref{eq:pdfp-step}}
  \Until{\eqref{eq:restart} holds \textbf{or} $T(z^{n,t})$ meets the termination criteria}
  \State $z^{n+1,0}\gets T(z^{n,t})$;\quad $n\gets n+1$
\Until{$z^{n,0}$ meets the termination criteria}
\State \Return $z^{n,0}$
\end{algorithmic}
\end{algorithm}

\subsection{Preconditioning}
\label{sec:pdfp-precond}

Like cuPDLPx, PDFP runs on diagonally preconditioned data, $\hat A=D_rAD_c$ and $x=D_c\hat x$. The column scaling changes the distance: $\varphi(x)=\sum_j q_j(D_c)_j|\hat x_j-\hat{\xt}_j|$. PDFP absorbs this change into the weights, $\hat q_j:=q_j(D_c)_j$, so the preconditioned problem has exactly the original distance as its objective. Only the box $Q$ changes, not the matrix $K$, so the step sizes are unaffected.

\subsection{Termination}
\label{sec:pdfp-term}

PDFP stops at the first PDHG point $z^{k+1}=T(z^k)$ whose relative primal residual, dual residual and gap fall below tolerances, as cuPDLPx does for LP:
\begin{align*}
\norm{\big((A_{\mathrm{ineq}}x^{k+1}-b_{\mathrm{ineq}})_+,\ A_{\mathrm{eq}}x^{k+1}-b_{\mathrm{eq}}\big)}_2&\le\varepsilon_r\,(1+\norm{b}_2) &&\text{(primal residual)},\\
\norm{s^{k+1}+A^\top y^{k+1}-r^{k+1}}_2&\le\varepsilon_g\,(1+\norm{q}_2) &&\text{(dual residual)},\\
\big|\varphi(x^{k+1})-L(x^{k+1},s^{k+1},y^{k+1})\big|&\le\varepsilon_g\,\big(1+\textstyle\sum_j q_j|\xt_j|\big) &&\text{(gap)}.
\end{align*}
Here $r^{k+1}:=\eta^{-1}(x^k-x^{k+1})+s^k+A^\top y^k\in N_X(x^{k+1})$ is the bound multiplier that the $x$-update recovers. If all three quantities vanish, $z^{k+1}$ is a saddle point of~\eqref{eq:pdfp}. The pump sets the tolerances (Section~\ref{sec:tolerance}).

%======================================================================
\section{Convergence of PDFP}
\label{sec:theory}

In this section, we show that PDFP inherits the accelerated linear convergence that r$^2$HPDHG attains on LP~\citep{HalpernPDHG2024}. That analysis rests on two properties of the PDHG operator $T$. First, $T$ is firmly nonexpansive in the $P$-norm, which holds for PDHG on every convex--concave saddle-point problem (Section~\ref{sec:pdhg}). Second, $T$ is \emph{sharp}: the fixed-point residual $\norm{z-T(z)}_P$ bounds the distance from $z$ to the solution set. \citet{HalpernPDHG2024} prove sharpness for LP from Hoffman's bound for the KKT system of the LP~\citep{Hoffman1952,ApplegateHinderLuLubin2023}. Problem~\eqref{eq:pdfp} is not the primal--dual form of an LP, because its dual set contains the box $Q$, so that proof does not apply directly. We prove sharpness for~\eqref{eq:pdfp} (Proposition~\ref{prop:sharp}), and the rate then follows as in \citet{HalpernPDHG2024} (Theorem~\ref{thm:rate}). Throughout, the target $\xt$, the weights $q$ and the step sizes are fixed, so the results describe one projection, and $T$ is the operator~\eqref{eq:pdfp-step}.

\begin{assumption}\label{ass:pdfp}
(i) $\Omega\neq\emptyset$. (ii) $\eta\sigma\norm{K}_2^2<1$.
\end{assumption}
Assumption~\ref{ass:pdfp}(i) holds whenever the pump runs, since the pump starts from a solution of the LP relaxation. Under Assumption~\ref{ass:pdfp}(ii), $P$ is positive definite. We denote by $Z^\star$ the set of saddle points of~\eqref{eq:pdfp}.

\begin{lemma}\label{lem:saddle}
Under Assumption~\ref{ass:pdfp}, $Z^\star$ is nonempty and $Z^\star=\Fix(T)$. A point $(x^\star,s^\star,y^\star)$ lies in $Z^\star$ if and only if $x^\star$ solves~\eqref{eq:proj} and $(s^\star,y^\star)$ solves the dual~\eqref{eq:pdfp-dual}.
\end{lemma}
\begin{proof}
Problem~\eqref{eq:proj} is a feasible linear program after the absolute values are linearized, and its objective is bounded below by $0$, so it has a solution. The functions in~\eqref{eq:pdfp} are polyhedral: $L$ is bilinear, and $X$, $Q$ and $Y$ are polyhedra. Polyhedral duality therefore gives a dual solution with zero duality gap, and the primal--dual solution pairs are exactly the saddle points~\citep[Cor.~31.2.1 and Sect.~36]{Rockafellar1970}. The saddle points are the zeros of $F$, which are the fixed points of $T$ by~\eqref{eq:ppm}.
\end{proof}

\begin{proposition}[sharpness]\label{prop:sharp}
Under Assumption~\ref{ass:pdfp}, for any $R>0$ there exists $\alpha>0$ such that
\[
\alpha\,\dist_P(z,Z^\star)\le\norm{z-T(z)}_P\qquad\text{for every $z$ with }\norm{T(z)}_2\le R.
\]
\end{proposition}
\begin{proof}
As in \citet[Prop.~2]{HalpernPDHG2024}, we first bound the distance to $Z^\star$ by the residual of the inclusion $0\in F(u)$, and then convert this bound to the fixed-point residual. Here $F(x,s,y)=\big(-s-A^\top y+N_X(x),\ x-\xt+N_Q(s),\ Ax-b+N_Y(y)\big)$ is the operator of~\eqref{eq:pdfp}.

\emph{Step 1.} The operator $F$ is polyhedral, i.e.\ its graph $\{(u,v):v\in F(u)\}$ is a finite union of polyhedra. Indeed, $X$, $Q$ and $Y$ are products of intervals, so $N_X$, $N_Q$ and $N_Y$ are products of one-dimensional normal-cone maps, each with a graph that is a finite union of polyhedra in $\R^2$. Adding the affine map $(x,s,y)\mapsto(-s-A^\top y,\ x-\xt,\ Ax-b)$ keeps the graph a finite union of polyhedra. Then $F^{-1}$ is polyhedral as well, and by \citet[Prop.~1]{Robinson1981} it is locally upper Lipschitz at $0$: there are $r_0>0$ and $\kappa_0$ such that $F^{-1}(v)\subseteq Z^\star+\kappa_0\norm{v}_2\mathbb{B}$ whenever $\norm{v}_2\le r_0$, where $\mathbb{B}$ is the unit ball. Now let $\norm{u}_2\le R$ and $v\in F(u)$, so that $u\in F^{-1}(v)$. If $\norm{v}_2\le r_0$, then $\dist_2(u,Z^\star)\le\kappa_0\norm{v}_2$. Otherwise $\dist_2(u,Z^\star)\le D_R\le(D_R/r_0)\norm{v}_2$, where $D_R:=\sup_{\norm{u'}_2\le R}\dist_2(u',Z^\star)$ is finite because $Z^\star\neq\emptyset$. With $\kappa_R:=\max\{\kappa_0,D_R/r_0\}$, taking the infimum over $v\in F(u)$ gives
\begin{equation}\label{eq:eb}
\dist_2(u,Z^\star)\le\kappa_R\,\dist_2\big(0,F(u)\big)\qquad\text{for every $u$ with }\norm{u}_2\le R.
\end{equation}

\emph{Step 2.} Let $u:=T(z)$ with $\norm{u}_2\le R$. By~\eqref{eq:ppm}, $P(z-u)\in F(u)$, and $\norm{Pw}_2\le\lambda_{\max}(P)^{1/2}\norm{w}_P$ for every $w$. With~\eqref{eq:eb},
\[
\dist_P(u,Z^\star)\le\lambda_{\max}(P)^{1/2}\dist_2(u,Z^\star)\le\kappa_R\lambda_{\max}(P)^{1/2}\norm{P(z-u)}_2\le\kappa_R\lambda_{\max}(P)\norm{z-u}_P.
\]
The triangle inequality then gives
\[
\dist_P(z,Z^\star)\le\norm{z-u}_P+\dist_P(u,Z^\star)\le\big(1+\kappa_R\lambda_{\max}(P)\big)\norm{z-T(z)}_P,
\]
so $\alpha:=\big(1+\kappa_R\lambda_{\max}(P)\big)^{-1}$ works.
\end{proof}

As for LP, the constant $\alpha$ can be very conservative. With Proposition~\ref{prop:sharp}, the accelerated rate of r$^2$HPDHG carries over to PDFP.

\begin{theorem}[accelerated linear convergence]\label{thm:rate}
Let Assumption~\ref{ass:pdfp} hold, let $\{z^{n,t}\}$ be generated by~\eqref{eq:halpern} with $T$ the operator~\eqref{eq:pdfp-step}, and restart with the fixed frequency $t_n=\lceil e/\alpha\rceil$, where $\alpha$ is the constant of Proposition~\ref{prop:sharp} for a radius $R$ with $\norm{T(z^{n,t})}_2\le R$ for all $n$ and $t$. Such an $R$ exists and depends only on $z^{0,0}$~\citep[Lemma~2]{HalpernPDHG2024}. Then
\[
\dist_P(z^{n,0},Z^\star)\le e^{-n}\dist_P(z^{0,0},Z^\star),\qquad n\ge0.
\]
Moreover, there exists $z^\infty\in Z^\star$ such that $\norm{z^{n,t}-z^\infty}_P\le2e^{-n}\dist_P(z^{0,0},Z^\star)$ for all $n\ge0$ and $0\le t\le t_n$, and $x^\infty$ solves~\eqref{eq:proj}.
\end{theorem}
\begin{proof}
The proof of \citet[Theorems~1 and~6]{HalpernPDHG2024} uses only two properties of $T$: firm nonexpansiveness in the $P$-norm, which also keeps the iterates bounded, and sharpness on a bounded set that contains them. With Proposition~\ref{prop:sharp} in place of their Proposition~2, it applies verbatim. The bound on $\norm{z^{n,t}-z^\infty}_P$ follows as in \citet[Theorem~1]{LuTang2026HPDLP}, and $x^\infty$ solves~\eqref{eq:proj} by Lemma~\ref{lem:saddle}.
\end{proof}

Theorem~\ref{thm:rate} implies that PDFP finds an $\varepsilon$-accurate solution within $O\big(\frac1\alpha\log\frac1\varepsilon\big)$ evaluations of $T$, the accelerated rate of r$^2$HPDHG for LP~\citep{HalpernPDHG2024}. The same two ingredients, firm nonexpansiveness and Proposition~\ref{prop:sharp}, carry the adaptive-restart analysis of \citet[Theorem~2]{HalpernPDHG2024} over to the restart condition~\eqref{eq:restart}; we omit the analogous proof.

%======================================================================
\section{Guarded two-phase constraint-propagation rounding}
\label{sec:cpr}

In this section, we present guarded two-phase CPR, the rounder of our pump. Sequential CPR (Algorithm~\ref{alg:seqcpr}) calls $\textsc{Propagate}$ after every variable it fixes by rounding, and with the projection on the GPU these calls take a large share of each pump iteration. On many instances most coordinates of $\xs_I$ are already integral, and sequential CPR usually leaves these \emph{LP-integral} coordinates at their LP values. Our rounder therefore groups the LP-integral coordinates into one \emph{batch} and calls $\textsc{Propagate}$ once for the whole batch. By default, it fixes these coordinates at their LP values; if earlier propagation has excluded an LP value, it uses the nearest value in the current domain. Only the fractional coordinates are rounded one at a time. When the batch succeeds, a call needs one propagation for the batch and at most one per fractional coordinate, instead of up to one per integer coordinate.

Fixing the batch at once can fail. Through the propagation of an earlier rounding, sequential CPR may move an LP-integral coordinate off its LP value, and a batch that fixes the coordinate at that value can then leave some row unsatisfiable (Example~\ref{ex:flip}). $\textsc{Propagate}$ detects such a row and raises its infeasibility flag (Section~\ref{sec:fp20}); we then say that the batch \emph{contradicts itself}. Our rounder checks the flag after the batch, a check we call the \emph{guard}. If the flag is raised, the rounder discards the batch and runs sequential CPR instead.

The rest of this section is organized as follows. Section~\ref{sec:cpr-alg} states the rounder. Section~\ref{sec:cpr-exact} proves that, in an idealized model of propagation, every rounding the rounder returns is one that sequential CPR returns for some order of the variables. Thus, under the stated assumptions, batching still returns a sequential CPR rounding, which need not be feasible. Sections~\ref{sec:cpr-schedule} and~\ref{sec:cpr-policy} describe how the pump chooses between two orders of processing the batch and the fractional coordinates, and when it tries the batch at all. Section~\ref{sec:l11} presents a rule that uses the constraints to choose the value of each rounding.

\begin{example}[fixing the batch at once can fail]\label{ex:flip}
Fixing the batch at once can fail even when sequential CPR finds a feasible rounding. Consider the integer program
\[
\min\ 0\qquad\text{s.t.}\quad a+2c\ge1,\qquad c-b\le\tfrac12,\qquad a,b,c\in\{0,1\}.
\]
Its LP relaxation replaces $a,b,c\in\{0,1\}$ by $a,b,c\in[0,1]$. The point $\xs=(a,b,c)=(0,0,\tfrac12)$ is optimal for the relaxation and satisfies both rows with equality. Its LP-integral coordinates are $a$ and $b$.

\emph{Sequential CPR.} The ranking places $a$ and $b$ first, since their fractionality is zero; take $a$ before $b$. CPR rounds $a$ to $0$ and propagates. The first row then forces $c\ge\tfrac12$, hence $c=1$, and the second row forces $b\ge\tfrac12$, hence $b=1$. Both $b$ and $c$ are now fixed by propagation and are read off their domains. The result $(0,1,1)$ is feasible, and propagation has moved $b$ off its LP value.

\emph{Batch.} The batch fixes $a=b=0$ at once and then propagates. The first row again forces $c=1$, and the second row then reads $1-0\le\tfrac12$, which cannot be satisfied. Together with $a=0$, the value $b=0$ is consistent only with a fractional $c$, so the batch contradicts itself. $\textsc{Propagate}$ raises its infeasibility flag on the second row, the guard discards the batch, and the rounder runs sequential CPR instead, which returns the feasible point $(0,1,1)$.
\end{example}

\subsection{The rounder}
\label{sec:cpr-alg}

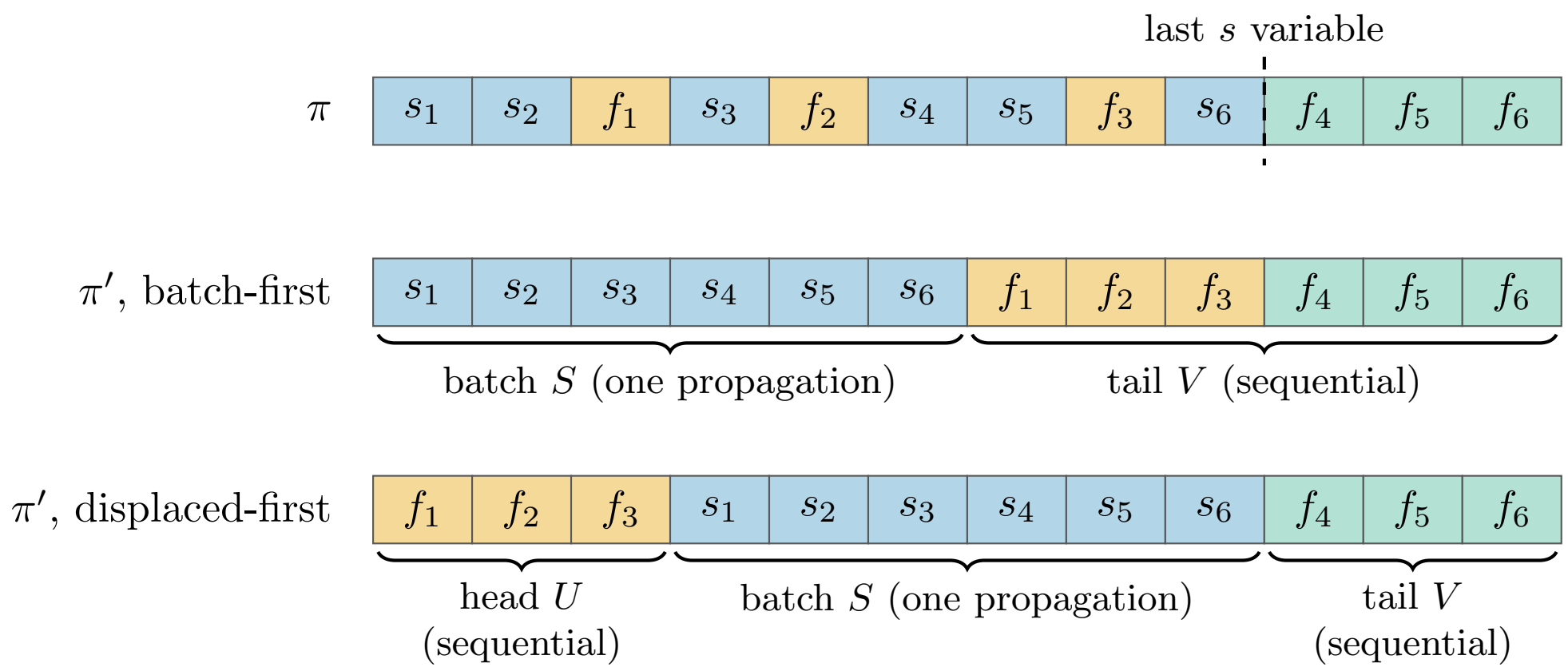
\begin{figure}[t]
\centering
\begin{tikzpicture}[x=0.82cm,y=1cm,
  cell/.style={draw=black!65, minimum width=0.82cm, minimum height=0.56cm, inner sep=0pt, font=\small},
  S/.style={cell, fill=schedS!30}, U/.style={cell, fill=schedU!40}, V/.style={cell, fill=schedV!30},
  lab/.style={anchor=east, font=\small, align=right},
  br/.style={decorate, decoration={brace, mirror, amplitude=4pt}, thick},
  bt/.style={anchor=north, font=\footnotesize, align=center}]
% the FP 2.0 ranking after the random swaps
\node[lab] at (-0.25,0) {$\pi$};
\foreach \v/\t [count=\i from 0] in {s_1/S,s_2/S,f_1/U,s_3/S,f_2/U,s_4/S,s_5/S,f_3/U,s_6/S,f_4/V,f_5/V,f_6/V}
  \node[\t] at (\i+0.5,0) {$\v$};
\draw[dashed, thick] (9,0.45) -- (9,-0.45);
\node[anchor=south, font=\footnotesize] at (9,0.45) {last $s$ variable};
% batch-first
\node[lab] at (-0.25,-1.5) {$\pi'$, batch-first};
\foreach \v/\t [count=\i from 0] in {s_1/S,s_2/S,s_3/S,s_4/S,s_5/S,s_6/S,f_1/U,f_2/U,f_3/U,f_4/V,f_5/V,f_6/V}
  \node[\t] at (\i+0.5,-1.5) {$\v$};
\draw[br] (0.05,-1.85) -- (5.95,-1.85) node[bt, midway, below=3pt] {batch $S$ (one propagation)};
\draw[br] (6.05,-1.85) -- (11.95,-1.85) node[bt, midway, below=3pt] {tail $V$ (sequential)};
% displaced-first
\node[lab] at (-0.25,-3.3) {$\pi'$, displaced-first};
\foreach \v/\t [count=\i from 0] in {f_1/U,f_2/U,f_3/U,s_1/S,s_2/S,s_3/S,s_4/S,s_5/S,s_6/S,f_4/V,f_5/V,f_6/V}
  \node[\t] at (\i+0.5,-3.3) {$\v$};
\draw[br] (0.05,-3.65) -- (2.95,-3.65) node[bt, midway, below=3pt] {head $U$\\(sequential)};
\draw[br] (3.05,-3.65) -- (8.95,-3.65) node[bt, midway, below=3pt] {batch $S$ (one propagation)};
\draw[br] (9.05,-3.65) -- (11.95,-3.65) node[bt, midway, below=3pt] {tail $V$\\(sequential)};
\end{tikzpicture}
\caption{The two schedules on an order $\pi$ in which some $f$ variables (orange) are interleaved among the $s$ variables (blue). The batch-first schedule places them after the $s$ variables, in the tail $V$. The displaced-first schedule places them before the $s$ variables, in the head $U$. The $f$ variables after the last $s$ variable (green) stay in the tail. Each block keeps its order in $\pi$.}
\label{fig:schedules}
\end{figure}

After the projection, sequential CPR sorts the integer variables by increasing fractionality of their values in $\xs$: the least fractional ones first, the most fractional ones last. We use $s$ to denote variables whose LP value is integral, and $f$ to denote variables whose LP value is fractional. The $s$ variables have fractionality zero, so they are ranked before the $f$ variables.\footnote{On instances with general integers, this holds among the binaries and among the general integers separately, because sequential CPR ranks all binaries first.} After the random swaps~\citep[Sect.~6]{FP20}, however, some $f$ variables may appear before the last $s$ variable. We call these the \emph{interleaved} fractional variables. We call this order after the swaps $\pi$. It is the order in which sequential CPR rounds the variables, and the first row of Figure~\ref{fig:schedules} shows an example.

Our rounder groups all $s$ variables into the batch $S$ and propagates their fixings together. It must therefore decide when to round the $f$ variables that are interleaved among them. We introduce two rounding schedules for this.

The \emph{batch-first} schedule takes all $f$ variables that are interleaved among the $s$ variables and places them after the block of $s$ variables. These variables, together with the untouched $f$ variables at the end, form a block which we call the \emph{tail} $V$. We then fix all variables in $S$ at their LP values and do one propagation. If the propagation results in a contradiction, we abandon the batch and run sequential CPR from the start, in the order $\pi$. Otherwise, we proceed to round the variables in $V$ one at a time, as sequential CPR does.

The \emph{displaced-first} schedule takes the same interleaved $f$ variables, but places them before the block of $s$ variables, where they form the \emph{head} $U$. The untouched $f$ variables at the end form the tail $V$. We first round the variables in $U$ one at a time. Their propagation narrows the domains of the $s$ variables, as it would in sequential CPR, and may already fix some of them. We then fix every $s$ variable that is not yet fixed and do one propagation. Each such variable is fixed at its LP value, or at the nearest value left in its domain if propagation has removed the LP value. As before, if the propagation results in a contradiction, we abandon the batch, together with the head, and run sequential CPR from the start. Otherwise, we round the variables in $V$ one at a time.

In both schedules, the rounder thus processes the variables in a new order, which we write $\pi'=U\cdot S\cdot V$: first $U$, then $S$, then $V$, where $U$ is empty for the batch-first schedule. Within each block, the variables keep their order in $\pi$. The check for a contradiction after the batch is the guard, and we call the case in which the guard does not fire the \emph{fast path}. On the fast path, a call needs one propagation for the batch and at most one for each $f$ variable. Algorithm~\ref{alg:cpr} summarizes the rounder, and Section~\ref{sec:cpr-schedule} describes how the pump chooses between the two schedules.

\begin{algorithm}[t]
\caption{Guarded two-phase CPR.}
\label{alg:cpr}
\begin{algorithmic}[1]
\Require LP point $\xs$; the order $\pi$ of sequential CPR; a schedule (batch-first or displaced-first)
\State $S\gets$ the $s$ variables, whose LP values are integral
\If{the schedule is displaced-first}
  \State $U\gets$ the $f$ variables that $\pi$ places before the last $s$ variable
\Else
  \State $U\gets\emptyset$ \Comment{batch-first}
\EndIf
\State $V\gets$ the $f$ variables not in $U$
\State round the variables in $U$ one at a time, in their order in $\pi$ \Comment{as in Algorithm~\ref{alg:seqcpr}}
\State fix each unfixed $s$ variable to its LP value if available, or to the nearest value in its current domain; propagate once
\If{the propagation results in a contradiction} \Comment{the guard}
  \State \Return the rounding of sequential CPR (Algorithm~\ref{alg:seqcpr}), run from the start in order $\pi$
\EndIf
\State round the variables in $V$ one at a time, in their order in $\pi$ \Comment{the fast path}
\State \Return the rounded vector $\xt$
\end{algorithmic}
\end{algorithm}

\subsection{Exactness of the guard}
\label{sec:cpr-exact}

The guard catches the failure of Example~\ref{ex:flip}. We now show that, in an idealized model of propagation, the guard is exact: the guard accepts the batch if and only if sequential CPR in the order $\pi'$ makes the same assignments on $S$ without emptying the box. When it accepts, both procedures also start the tail from the same box and therefore return the same rounding.

\paragraph{The model.} A \emph{box} $B=\prod_jB_j$ assigns each variable an interval, its domain; the domain of an integer variable is an integer interval. Each row $i$ defines a \emph{propagator} $p_i$, the bound strengthening of row~$i$ (Section~\ref{sec:fp20}), which maps a box $B$ to a box $p_i(B)\subseteq B$. Propagators are monotone: $B\subseteq B'$ implies $p_i(B)\subseteq p_i(B')$~\citep[Sect.~3.2]{FP20}. A box that no propagator changes is a \emph{fixed point}; here a fixed point is a box, not a single point. We compare boxes by inclusion, and let $\Prop(B)$ be the largest fixed point inside $B$: the fixed point inside $B$ that contains every other fixed point inside $B$, possibly the empty box. It exists by Tarski's fixed-point theorem~\citep{Tarski1955}, and $\Prop$ is monotone as well. With finite domains, applying the propagators in any order until none changes the box ends at $\Prop(B)$~\citep{Apt2003}. So $\Prop(B)$ is the result of propagating from $B$, whatever the order in which the rows are processed. In this model, sequential CPR in an order $(j_1,\dots,j_N)$ from a box $A_0$ chooses, for $k=1,\dots,N$, a value $v_k$ from the current box $A_{k-1}$ and sets $A_k:=\Prop(A_{k-1}\cap\{x_{j_k}=v_k\})$. The \emph{plain} rule of FP~2.0 chooses $v_k$ as the value of a singleton domain, and otherwise as the integer in the domain nearest to $\xs_{j_k}$ (lines 5 and~7 of Algorithm~\ref{alg:seqcpr}).

\begin{assumption}\label{ass:cpr}
(i) $\textsc{Propagate}$ computes $\Prop$ exactly, and raises its infeasibility flag exactly when $\Prop$ is empty. (ii) Arithmetic is exact. (iii) The initial box $B_0$ set in line~3 of Algorithm~\ref{alg:seqcpr} is a fixed point, $\Prop(B_0)=B_0$.
\end{assumption}

To state the result, let $A_U$ be the box after the head, $S_U\subseteq S$ the members of $S$ that are not yet fixed in $A_U$, and $t_j$ the plain value of $j\in S_U$ in $A_U$, which is $\xs_j$ unless the head has excluded it. Fixing the batch and propagating gives the box
\[
C:=\Prop\big(A_U\cap\{x:\ x_j=t_j\ \text{for all } j\in S_U\}\big).
\]
In the model, the batch contradicts itself exactly when $C$ is empty.

\begin{theorem}[the guard is exact]\label{thm:cpr}
Let Assumption~\ref{ass:cpr} hold. Use the plain rule on $S$, and on $U\cup V$ any rule that chooses each value from the current box alone. Then $C\neq\emptyset$ if and only if sequential CPR in the order $\pi'=U\cdot S\cdot V$, processing $S$ from $A_U$, assigns every $j\in S_U$ the value $t_j$ without emptying the box. In that case, Algorithm~\ref{alg:cpr} returns the same rounding as sequential CPR in the order $\pi'$, for any order of $S$ within $\pi'$.
\end{theorem}
\begin{proof}
Both runs process $U$ in the same way and reach $A_U$, which is a fixed point (by (iii) if $U=\emptyset$). Let the sequential run process $S$ in the order $s_1,\dots,s_N$, with values $v_1,\dots,v_N$ and boxes $A'_0=A_U,A'_1,\dots,A'_N$.

\emph{Only if, and equal roundings.} Let $C\neq\emptyset$. We show by induction on $k$ that $C\subseteq A'_k$ and that $v_k$ is the value that the batch gives $s_k$. The base case $C\subseteq A'_0$ holds because $\Prop(B)\subseteq B$. Assume $C\subseteq A'_{k-1}$. If $s_k\notin S_U$, its domain is still the singleton it had in $A_U$, both runs read that value, and $A'_k=A'_{k-1}$. If $s_k\in S_U$, every point of the nonempty box $C$ has $x_{s_k}=t_{s_k}$, so $t_{s_k}$ lies in the domain $D$ of $s_k$ in $A'_{k-1}$. Since $D$ is contained in the domain of $s_k$ in $A_U$, and $t_{s_k}$ is the integer of that larger domain nearest to $\xs_{s_k}$, it is also the integer of $D$ nearest to $\xs_{s_k}$, under the same tie-breaking rule. Hence $v_k=t_{s_k}$, so $C\subseteq A'_{k-1}\cap\{x_{s_k}=v_k\}$, and the monotonicity of $\Prop$ gives
\[
C=\Prop(C)\subseteq\Prop\big(A'_{k-1}\cap\{x_{s_k}=v_k\}\big)=A'_k.
\]
So every $j\in S_U$ receives $t_j$, and no box $A'_k$ is empty, since each contains $C$. After the last step, $A'_N$ is a fixed point inside $A_U\cap\{x:x_j=t_j\ \text{for all } j\in S_U\}$, so $A'_N\subseteq C$, because $C$ is the largest such fixed point. Hence $A'_N=C$. The tail $V$ then starts from the same box, in the same order and under the same rule, so the two runs return the same rounding.

\emph{If.} Suppose the sequential run assigns every $j\in S_U$ the value $t_j$ without emptying the box. Then $A'_N$ is a nonempty fixed point inside $A_U\cap\{x:x_j=t_j\ \text{for all } j\in S_U\}$, so $C\supseteq A'_N\neq\emptyset$.
\end{proof}

Thus, on the fast path, batching reproduces sequential CPR in the reordered sequence $\pi'$. In the model, the guard rejects a batch precisely when sequential CPR cannot make all of its assignments without emptying the box, as in Example~\ref{ex:flip}.

In the model, every rounding that Algorithm~\ref{alg:cpr} returns is therefore a sequential CPR rounding: in the order $\pi'$ on the fast path, and in the order of Algorithm~\ref{alg:seqcpr} when the guard fires. The two orders coincide when $\pi$ ranks all LP-integral variables first. This holds on pure-binary instances before the random swaps start, and there every call returns FP~2.0's own rounding. Once the swaps mix LP-integral and fractional variables, the fast path may return a different sequential CPR rounding.

\subsection{Choosing the schedule}
\label{sec:cpr-schedule}

The two schedules differ only when $\pi$ ranks some fractional variable before an LP-integral one. The displaced-first schedule can help when the LP objective does not strongly favor keeping the LP-integral coordinates at their values. Fixing these coordinates together may lead to a propagation contradiction, and rounding the head first lets propagation adjust them before the batch is fixed. We call $j\in S$ \emph{pushed} if its reduced cost in the root relaxation, the LP that the pump solves first, is nonzero with the sign that holds $x_j$ at its bound, so that moving $x_j$ off $\xs_j$ costs objective value. Once per run, the pump chooses displaced-first if the fraction of pushed coordinates in $S$ is below a threshold $\theta_{\rm push}$, and batch-first otherwise. In our experiments, displaced-first led to more solves only in the first case.

\subsection{Learning, per run, whether to batch}
\label{sec:cpr-policy}

A batch that contradicts itself wastes work: the call pays for the batch and then for a full run of sequential CPR. How often batches contradict themselves depends on the instance and changes during a run. Before each rounding, the pump therefore chooses between \emph{trying} the batch, that is, calling Algorithm~\ref{alg:cpr}, and \emph{skipping} it, that is, calling sequential CPR directly. It maintains three exponentially weighted mean costs observed so far in the run: $\bar s$ for sequential CPR, $\bar h$ for a call whose batch is accepted, and $\bar m$ for the work done before falling back after a rejected batch. Each mean gives weight $\xi$ to its newest observation. It also keeps the fraction $p$ of tried batches that contradicted themselves, discounted by a factor $\nu$ at every tried call. It tries the batch if and only if
\[
p\,(\bar m+\bar s)+(1-p)\,\bar h<\bar s.
\]
The left-hand side estimates the cost of trying: a rejected batch costs $\bar m$ plus the fallback cost $\bar s$, whereas an accepted batch costs $\bar h$. The right-hand side estimates the cost of skipping. Since only a tried batch shows whether it contradicts itself, the rule tries the batch again after $R_p$ consecutive skips.

\subsection{Violation-guided rounding}
\label{sec:l11}

The plain rule rounds $\xs_j$ to the nearest integer in the current domain: CPR uses the constraints to narrow each domain, but not to choose a value inside it. Our \emph{violation-guided rule} instead sets $x_j$ to the integer in its current domain $[\underline v,\overline v]$ that minimizes the change in constraint violation computed from the variables fixed so far. Let $D$ be the set of integer variables whose domain in the current box is a singleton, that is, every integer variable fixed so far in the call, whether by a rounding, by the batch or by propagation. Let $\mathrm{act}_i:=\sum_{k\in D}a_{ik}\xt_k$ be the \emph{running activity} of row $i$, the part of its left-hand side that these variables contribute. The rule chooses
\[
v_j\in\argmin_{v\in[\underline v,\overline v]\cap\Z}\ \mathrm{score}(v),\qquad
\mathrm{score}(v):=\sum_{i:\,a_{ij}\neq0}\Big[\mathrm{viol}_i\big(\mathrm{act}_i+a_{ij}v\big)-\mathrm{viol}_i(\mathrm{act}_i)\Big],
\]
where $\mathrm{viol}_i(r):=(r-b_i)_+$ for an inequality row and $|r-b_i|$ for an equation. So $\mathrm{score}(v)$ is the increase in the violation of the rows that contain $x_j$ when $x_j=v$ is added to the running activity. Ties are broken toward $\xs_j$, so the rule keeps the plain value whenever the plain value is a minimizer.

Because $D$ contains every fixed variable, the rule depends only on the current box, so Theorem~\ref{thm:cpr} holds with the rule active. The rule acts only where variables are rounded one at a time, that is, in the head, the tail and sequential CPR, and never on the batch.

Used on every run, the rule helps on some instances and hurts on others. The pump therefore chooses once per run whether to use it, from the rounding of the root relaxation. It computes that rounding without the rule, and enables the rule for the rest of the run if the total violation $\sum_i\mathrm{viol}_i(A_i\xt)$ of that rounding reaches a threshold $\theta_{\rm viol}$. In our experiments, total violation was a better indicator of whether the rule would help than maximum violation.

%======================================================================
\section{Further enhancements to the GPU pump}
\label{sec:enhancements}

In this section, we describe further enhancements that control how the pump runs the projector and the rounder. Section~\ref{sec:gpu-state} keeps the projector and its data on the GPU for the whole run and warm-starts every projection. Section~\ref{sec:tolerance} sets the projector's tolerance, which trades the accuracy of each projection against the number of pump iterations. Section~\ref{sec:early-exit} stops the pump as soon as a PDHG iterate rounds to a feasible point, without finishing the projection. Section~\ref{sec:hashing} implements the pump's cycle test with less memory and less time. Section~\ref{sec:racing} runs two trajectories of the pump concurrently on one GPU, in the device time left idle while the pump rounds.

\subsection{A GPU-resident projector}
\label{sec:gpu-state}

The pump alternates a GPU phase, the projection, and a CPU phase, the rounding. We keep the constraint data and preconditioning on the GPU for the whole run, and retain the projector's iterates between projections. Each projection therefore requires only copying the target to the GPU and the projected point back to the CPU. Each projection is warm-started from the last iterate of the previous one, since between consecutive projections only the target changes.

\subsection{The tolerance policy}
\label{sec:tolerance}

Because a first-order solve becomes more expensive as the requested accuracy increases, the projector's tolerance largely determines the speed of a GPU pump. A fixed tight tolerance can spend too much time refining projections beyond the accuracy needed by the next rounding step. A fixed loose tolerance can return points too far from the optimal face, increasing the number of pump iterations. We use the geometric schedule of Scylla~\citep{Scylla2023}: the tolerance starts loose, shrinks by a constant factor at every projection, and stops at a floor. For PDFP the schedule sets $\varepsilon_g$, and $\varepsilon_r=\rho\,\varepsilon_g$ for a fixed ratio $\rho$. Table~\ref{tab:params} lists the parameters of our pump.

\begin{table}[t]
\centering
\caption{Parameters of our pump. Values taken from prior work name their source; the others are ours. Where the configuration with PDLP (Section~\ref{sec:configs}) differs, its value is in parentheses.}
\label{tab:params}
\begin{tabular}{@{}lll@{}}
\toprule
component & parameter & value \\
\midrule
FP loop & flip parameter $T$~\citep{FGL2005} & $20$ \\
 & restart period $R$~\citep{FGL2005} & $100$ \\
 & priority swap (FP~2.0) & $\lfloor|I|/10\rfloor$ transpositions from iteration $10$ \\
PDFP & step sizes $\eta=\lambda/\omega$, $\sigma=\lambda\omega$ & $\lambda=0.998/\sqrt{1+\norm{A}_2^2}$, $\omega$ by PID \\
 & restart triggers & cuPDLPx defaults \\
Tolerances & gap schedule: start, factor (Scylla) & $10^{-2}$, $0.98$ \\
 & gap floor (ours) & $10^{-4}$ ($10^{-6}$) \\
 & primal-residual ratio $\rho$ & $1/10$ ($1$) \\
 & root relaxation, solved by PDLP & $10^{-2}$ \\
Guarded CPR & schedule threshold $\theta_{\rm push}$ & $1\%$ pushed \\
 & policy: $\xi$, $\nu$, $R_p$ & $0.3$, $0.9$, $20$ (not used) \\
\begin{tabular}[t]{@{}l@{}}Violation-guided\\rounding\end{tabular} & gate $\theta_{\rm viol}$ (total violation) & $5000$ \\
\bottomrule
\end{tabular}
\end{table}

\subsection{Stopping at an integer-feasible iterate}
\label{sec:early-exit}

Whenever the projector checks its termination criteria, the pump also rounds the current iterate to the nearest integers and tests whether the result is feasible. If it is, the pump returns that point and stops without finishing the projection. If the test finds no feasible point, the projection continues unchanged.

\subsection{Cycle detection by hashing}
\label{sec:hashing}

The pump remembers every visited integer point to detect cycles (Algorithm~\ref{alg:fp}). The reference implementation~\citep{FeasPumpCollection} stores the points themselves and compares a candidate with each of them, so its memory and the cost of a test grow with the number of iterations. We store only a 64-bit hash of each visited point in a hash set. Each test requires one pass over the integer coordinates to compute the hash, followed by a hash-set lookup with expected constant cost. A hash collision can only trigger an unnecessary perturbation; it never causes an infeasible point to be accepted.

\subsection{Concurrent trajectories on one GPU}
\label{sec:racing}

Because our rounder runs on the CPU, the GPU idles during every rounding. The rounding is randomized, and on some instances the iteration count varies by orders of magnitude across seeds. Concurrent trajectories exploit this variation. We start two pump processes with independent seeds on the same GPU, take the first solution, and terminate the other process. The projections of each process tend to fall in the rounding phase of the other, so the two share the device with little contention when the rounding dominates. With little contention, the time to a solution is approximately the minimum of the two independent run times. This can substantially reduce the average time on instances with large variation across seeds. Concurrent trajectories thus let the pump terminate sooner on average (Section~\ref{sec:results}).

%\subsection{Mixed-integer instances}
%\label{sec:mixed}
%\todo{Describe what changes on models with continuous variables, once that configuration is frozen. Currently: when a rounding is integral but infeasible only through the continuous variables, the projection is re-solved with a tenfold tighter tolerance, up to five times per run, as in the reference code; PDFP computes its primal residual in original units; the propagation budget is off; the gate of the violation-guided rule counts the violation of the rows without continuous variables, or of all rows on models with general integers; and the two-stage scheme of Bertacco et al. runs where general integers outnumber binaries, while pure-integer instances run the single-stage loop of Algorithm~\ref{alg:fp}. Then update Sections~\ref{sec:pdfp-term} and~\ref{sec:l11} and Table~\ref{tab:params}, which describe the pure-integer configuration.}

%======================================================================
\section{Computational results}
\label{sec:results}

In this section, we evaluate our pump by the time it takes to find a first feasible solution on MIPLIB~2017 instances. Section~\ref{sec:setup} describes the benchmark and the protocol, Section~\ref{sec:baseline} the baselines, Section~\ref{sec:configs} the configurations of our pump, and Section~\ref{sec:results-main} the results.

\subsection{Setup}
\label{sec:setup}

\paragraph{Benchmark dataset.} We conduct experiments on the large instances of the MIPLIB~2017 collection~\citep{MIPLIB2017}. We focus on large instances because GPU implementations of PDHG are designed for large-scale problems~\citep{cuPDLP2023,cuPDLPx2025}. The collection has 159 instances with at least $10^6$ nonzeros in the constraint matrix. We exclude 4 of them with indicator constraints, which are not of the form~\eqref{eq:milp}, 6 that are infeasible, and 2 whose feasibility is unknown. Of the remaining 147 instances, we exclude 5 for which we could not obtain a verified presolved model (see below). This leaves 142 instances, 53 pure-integer and 89 mixed-integer after presolve. Table~\ref{tab:benchmark} groups them by size.

\begin{table}[htb]
\centering
\begin{tabular}{lrrr}
\toprule
Number of nonzeros & 1M--10M & $>$10M & Total \\
\midrule
Pure integer & 43 & 10 & 53 \\
Mixed integer & 78 & 11 & 89 \\
\midrule
Total & 121 & 21 & 142 \\
\bottomrule
\end{tabular}
\caption{Size and type of the 142 benchmark instances. The number of nonzeros is that of the original constraint matrix; the type is that of the presolved model.}
\label{tab:benchmark}
\end{table}

\paragraph{Presolve.} We presolve every instance once, before the experiments, with SCIP~10.0.3~\citep{SCIP}, using the MIPLIB~2017 trivial-presolve settings file~\citep{MIPLIB2017}. We turn off symmetry detection and dual reductions, since both can cut off feasible solutions. For each presolved model, we check that its solutions map back to solutions of the original instance, and that it admits a known feasible solution of the original instance. All methods read the same presolved model, and the time of presolve is not counted.

\paragraph{Software.} We implement our pump in Julia~1.12.6~\citep{Bezanson2017Julia}, with CUDA.jl~6.1.0~\citep{Besard2019CUDAjl} as the interface to NVIDIA GPUs. The projections run on the GPU, and the rest of the pump, including the rounder, runs on one CPU thread. To exclude compilation time, each process first runs the pump on a small instance before timing the benchmark instance.

\paragraph{Computing environment.} We use one server with NVIDIA H100 80GB HBM3 GPUs, NVIDIA driver 595.71.05 with CUDA~13.2, two Intel Xeon Platinum 8481C CPUs with 52 cores each, and 1.8~TB of memory. Each run has exclusive use of one GPU, shared only by its two processes when it runs two concurrent trajectories.

\paragraph{Seeds and concurrent trajectories.} The pump is randomized, and its time varies with the seed. We therefore run each instance with four seeds, which gives 568 runs per method. A run with concurrent trajectories (Section~\ref{sec:racing}) starts two processes with independent seeds on its GPU and stops when either finds a solution.

\paragraph{Feasibility check.} A run solves its instance if, within the time limit, it returns a point that satisfies integrality, the bounds, and every constraint of the presolved model to an absolute tolerance of $10^{-6}$. Because the methods use different internal feasibility tolerances, we apply this common check to every returned point.

\paragraph{Time limit.} We apply a time limit of 1800 seconds to every run. For our pump and for the FP~2.0 reference implementation, the time starts when the pump is called on the loaded model, so it includes the solve of the root LP relaxation and excludes reading the model.

\paragraph{Shifted geometric mean.} We use the shifted geometric mean of the time to a first feasible solution to evaluate performance across a collection of runs. Formally, the shifted geometric mean is $\big(\prod_{i=1}^n (t_i+\Delta)\big)^{1/n}-\Delta$, where $t_i$ is the time of the $i$-th of $n$ runs. We shift by $\Delta=10$ seconds and denote this metric by SGM10. If a run does not solve its instance, its time is set to the time limit. We report the number of solved runs and SGM10 for the pure-integer and mixed-integer subsets and for all runs combined.

\subsection{Baselines}
\label{sec:baseline}

We compare our pump with two external methods. The first is the FP~2.0 reference implementation of Salvagnin and Mexi~\citep{FeasPumpCollection} (commit 79515d6), with its LP solver replaced by cuPDLPx~0.3.0 (commit 8b43b8c)~\citep{cuPDLPx2025}, kept on the GPU across projections. The second is NVIDIA's cuOpt~25.10.1~\citep{Corduk2025cuOpt}, run with its full primal-heuristic portfolio and timed from the start of its solve to its first incumbent.

\subsection{Configurations of our pump}
\label{sec:configs}

We report two configurations of our pump. We run each either with a single process or with two concurrent trajectories (Section~\ref{sec:racing}), giving the four rows for our pump in Table~\ref{tab:results}. Both configurations use the guarded two-phase rounder with its choice of schedule and the violation-guided rule (Section~\ref{sec:cpr}). Both keep the projector on the GPU (Section~\ref{sec:gpu-state}), follow the tolerance schedule of Section~\ref{sec:tolerance}, and detect cycles by hashing (Section~\ref{sec:hashing}). They differ as follows.

\paragraph{Enhanced FP (PDFP).} This configuration runs the single-stage loop of Algorithm~\ref{alg:fp} with PDFP as its projector (Section~\ref{sec:pdfp}). Its rounder learns, per run, whether to try the batch (Section~\ref{sec:cpr-policy}). It does not use the mid-projection feasibility check of Section~\ref{sec:early-exit}.

\paragraph{Enhanced FP (PDLP).} This configuration is an earlier version of our pump. Its projector is our PDLP projector, which solves the LP form~\eqref{eq:proj-lp} of the projection. Its rounder tries the batch at every rounding. It uses the mid-projection feasibility check of Section~\ref{sec:early-exit}. On instances with general integers, it runs the first two stages of the scheme of \citet{BFL2007}: it pumps first on the binaries and then on all integer variables. It also differs from the PDFP configuration in some parameters and in further details of the rounder.

\subsection{Results}
\label{sec:results-main}

Table~\ref{tab:results} reports the number of solved runs and SGM10 of each method on the pure-integer instances, on the mixed-integer instances, and on all instances.

\begin{table}[htb]
\centering
\setlength{\tabcolsep}{5pt}
\begin{tabular}{@{}lcccccc@{}}
\toprule
 & \multicolumn{2}{c}{Pure integer (212 runs)} & \multicolumn{2}{c}{Mixed integer (356 runs)} & \multicolumn{2}{c}{All (568 runs)} \\
\cmidrule(lr){2-3}\cmidrule(lr){4-5}\cmidrule(l){6-7}
Method & Solved & SGM10 & Solved & SGM10 & Solved & SGM10 \\
\midrule
FP~2.0 with cuPDLPx & 179 & 48.8 & 202 & 234.2 & 381 & 133.6 \\
NVIDIA cuOpt & 190 & 73.2 & 211 & 340.9 & 401 & 195.0 \\
\midrule
Enhanced FP (PDLP) & 190 & 39.9 & 190 & 341.7 & 380 & 159.7 \\
\quad + concurrent trajectories & 192 & 36.3 & 200 & 306.7 & 392 & 144.5 \\
Enhanced FP (PDFP) & 193 & 32.2 & 226 & 236.3 & 419 & 117.5 \\
\quad + concurrent trajectories & \textbf{194} & \textbf{31.3} & \textbf{230} & \textbf{226.6} & \textbf{424} & \textbf{113.4} \\
\bottomrule
\end{tabular}
\caption{Solved runs and SGM10 in seconds on the 142 benchmark instances, with four seeds each. A run with concurrent trajectories starts two processes with independent seeds on one GPU (Section~\ref{sec:racing}); the other runs use one process. The best value in each column is in bold.}
\label{tab:results}
\end{table}

With concurrent trajectories, our pump with PDFP solves the most runs and has the lowest SGM10 in each instance group of Table~\ref{tab:results}. On the pure-integer instances, every method solves most runs, and the main gain is in time: the PDFP configuration with concurrent trajectories has an SGM10 of 31.3 seconds, against 48.8 for FP~2.0 with cuPDLPx and 73.2 for cuOpt. On the mixed-integer instances, the main gain is in the number of solved runs: 230, against 202 and 211. With one process per run, our PDFP configuration still solves the most runs in each instance group. It also has the lowest SGM10 except on the mixed-integer instances, where FP~2.0 with cuPDLPx has a slightly lower SGM10 (234.2 against 236.3 seconds).

Running two concurrent trajectories adds 5 solved runs to the configuration with PDFP and 12 to the one with PDLP, and lowers their SGM10 from 117.5 to 113.4 seconds and from 159.7 to 144.5 seconds. With one process per run, the configuration with PDFP solves 39 more runs than the one with PDLP, 36 of them on the mixed-integer instances.

\section{Conclusion}
\label{sec:conclusion}

In this paper, we present an end-to-end primal heuristic based on Feasibility Pump and PDHG that runs on GPU. Its projector, PDFP, is a first-order method that solves the pump's projection problem directly, without auxiliary primal variables or constraints, and we prove that it keeps the accelerated linear rate of restarted Halpern PDHG. Its rounder, guarded two-phase CPR, introduces batch fixing and violation-guided rounding decisions on top of the original CPR. Running two concurrent trajectories on one GPU further speeds up the pump. Early computational experiments on large MIPLIB~2017 instances show that our pump solves more runs than both the FP~2.0 reference implementation with cuPDLPx and NVIDIA's cuOpt. In shifted geometric mean, it finds a first feasible solution 1.18$\times$ faster than the FP~2.0 reference implementation and 1.72$\times$ faster than cuOpt. These results highlight the effectiveness of these algorithmic enhancements to the feasibility pump in finding the initial feasible solution.

\section*{Acknowledgments}
Ruofeng Liu is supported by the National Science Foundation Graduate Research Fellowship Program under Grant No.~2141064. Any opinions, findings, and conclusions or recommendations expressed in this material are those of the authors and do not necessarily reflect the views of the National Science Foundation. Haihao Lu is partially supported by AFOSR Grant No. FA9550-24-1-0051, ONR Grant No. N000142412735, and a Sloan Research Fellowship. Rahul Mazumder acknowledges research support from ONR N000142512504
and U.S. Department of Energy, Office of Science, Award Number DE-SC0026379.
We gratefully acknowledge NVIDIA Academic Grant Program for providing GPU computing resources.
\bibliographystyle{plainnat}
\bibliography{references}

\end{document}